\documentclass[12pt,reqno]{article}

\usepackage[usenames]{color}
\usepackage{amssymb}
\usepackage{graphicx}
\usepackage{amscd}

\usepackage[colorlinks=true,
linkcolor=webgreen,
filecolor=webbrown,
citecolor=webgreen]{hyperref}

\definecolor{webgreen}{rgb}{0,.5,0}
\definecolor{webbrown}{rgb}{.6,0,0}

\usepackage{color}
\usepackage{fullpage}
\usepackage{float}

\usepackage{graphics,amsmath,amssymb}
\usepackage{amsthm}
\usepackage{amsfonts}
\usepackage{latexsym}
\usepackage{epsf}

\begin{document}
	
	\theoremstyle{plain}
	\newtheorem{theorem}{Theorem}
	\newtheorem{corollary}[theorem]{Corollary}
	\newtheorem{lemma}[theorem]{Lemma}
	\newtheorem{proposition}[theorem]{Proposition}
	\newtheorem{question}[theorem]{Question}
	
	\theoremstyle{definition}
	\newtheorem{definition}[theorem]{Definition}
	\newtheorem{example}[theorem]{Example}
	\newtheorem{conjecture}[theorem]{Conjecture}
	
	\theoremstyle{remark}
	\newtheorem{remark}[theorem]{Remark}
	
	\begin{center}
		\vskip 1cm{\LARGE\bf Odd Spoof Multiperfect Numbers}
		\vskip 1cm
		\large
		L\'aszl\'o T\'oth \\
		Eischen, Grand Duchy of Luxembourg \\
		\href{mailto:uk.laszlo.toth@gmail.com}{\tt uk.laszlo.toth@gmail.com}
	\end{center}
	
	\vskip .2 in
	
	\begin{abstract}
		We generalize the definition of spoof perfect numbers to multiperfect numbers and study their characteristics. As a result, we find several new odd spoof multiperfect numbers, akin to Descartes' number. An example is $8999757$, which would be an odd multiperfect number, if one of its prime factors, $61$, was a square. We briefly describe an algorithm for searching for such numbers and discuss a few of their properties. 
	\end{abstract}

\section{Introduction}

Let $\sigma(n)$ denote the sum of divisors of the positive integer $n$. If $\sigma(n) = 2n$, $n$ is said to be \textit{perfect}. While many examples of even perfect numbers are known, no odd perfect numbers have been found so far. Their existence is the subject of extensive research, and the abundance of constraints that an an odd perfect number must satisfy has left many skeptical about their existence.

Descartes however noted in 1638 that the number
$$
\mathcal{D}=198585576189
$$
would be an odd perfect number if only one of its composite factors were prime. Indeed, $\mathcal{D}=3^2 \cdot 7^2 \cdot 11^2 \cdot 13^2 \cdot 22021$ and assuming that $22021$ is prime, we have
\begin{align*}
\sigma(\mathcal{D}) &= (3^2+3+1) \cdot (7^2+7+1)\cdot(11^2+11+1)\cdot(13^2 +13+1)\cdot(22021+1) \\ &=2\mathcal{D}.
\end{align*}
Alas, $22021=19^2 \cdot 61$ so $\mathcal{D}$ is not perfect. However, Descartes' example motivated the search for so-called ``spoof perfect'' numbers, integers $s=nx$ with the property
$$
\frac{\sigma(n)}{2n} = \frac{x}{x+1}
$$
for positive integers $n$ and $x$ such that $n,x\geq2$. If $x$, the so-called \textit{spoof factor}, is a prime number, then $nx$ is perfect. On the other hand, if both $n$ and $x$ are odd, but $x$ is not prime, then we have a \textit{Descartes number}, also referred to as an \textit{odd spoof perfect number}. Descartes' example is the only currently known odd spoof perfect number.

Considerable research has been carried out in finding examples similar to $\mathcal{D}$. For instance, T\'oth \cite{Toth21} showed that if $s=nx$ denotes a spoof perfect number with pseudo-prime $x$, other than Descartes' example, then $n$ must be greater than $10^{12}$. Others, like Voight \cite{Voight}, generalized the concept of spoof perfect numbers by allowing negative spoof factors. One such example, due to Voight, is
$$
\mathcal{V} = 3^4 \cdot 7^2 \cdot 11^2 \cdot 19^2 \cdot (-127),
$$
which would be perfect if only $-127$ were positive. The BYU Computational Number Theory Group \cite{Pace spoof group} further relaxed these conditions. They discovered that the number
$$
11025 = 1 \cdot 9 \cdot 25 \cdot 49 = (1)^2 \cdot (-3)^2 \cdot (-5)^2 \cdot (49)
$$
would be an odd perfect number if $1$, $-3$, $-5$, and $49$ were all assumed to be prime.

\section{Scope of This Paper}

In this paper we provide a generalization of spoof perfect numbers to multiperfect (or ``multiply perfect'') numbers, integers $n$ such that $\sigma(n)=kn$ for some positive integer $k>2$. Indeed, motivated by the fact that no odd spoof perfect number $s=nx$ with the property
$$
\frac{\sigma(n)}{2n} = \frac{x}{x+1}
$$
has been found other than Descartes' example, we turn our attention instead to numbers such that
$$
\frac{\sigma(n)}{kn} = \frac{x}{x+1},
$$
i.e., the ``spoof'' equivalent of multiperfect numbers, and
$$
\frac{\sigma(n)}{kn} = \frac{x}{x^2+x+1},
$$
numbers that would be multiperfect if only one of their factors $x$ was the square of a prime. Our methods allow us to find several new odd spoof multiperfect numbers. For instance, let 
$$
S = 8999757 = 3^2 \cdot 13^2 \cdot 61 \cdot 97.
$$
If one assumes (wrongly) that $61$ is a square, we have :
\begin{align*}
\sigma(S) &= (3^2 + 3 + 1) (13^2 + 13 + 1) (97 + 1) (61^2 + 61 + 1) \\
&= (13) \cdot (3\cdot61) \cdot (2 \cdot 7^2) \cdot ( 3 \cdot 13 \cdot 97) \\
&= 2 \cdot 3^2 \cdot 7^2 \cdot 13^2 \cdot 61 \cdot 97 \\
&= 98 \cdot 3^2 \cdot 13^2 \cdot 61 \cdot 97 \\
&= 98 S.
\end{align*}
And thus, $S$ would be an odd $98$-perfect number if only its spoof factor ($61$) was a square. 

Inspired by this result, we conducted a computer search in hopes of finding other examples of this kind. And indeed, we found several other odd spoof multiperfect numbers, which we will present further on in this paper. In the next sections, we will define spoof $k$-perfect numbers in a more rigorous manner and provide a few examples, both odd and even. We will then describe our algorithm that allows to quickly find such numbers and conclude with a few remarks about Robin's inequality and spoof multiperfect numbers.

\section{Spoof Multiperfect Numbers}

We begin this section by defining spoof $k$-perfect numbers as follows. 
\begin{definition}[Spoof $k$-perfect number] \label{def1}
Let $s=nx$ be a positive integer such that $n,x \in \mathbb{N}$ and $n,x\geq2$. Then $s$ 
\begin{enumerate}
    \item is a \textit{spoof $k$-perfect number of the first kind} if $\sigma(n) (x+1) = knx$,
    \item is a \textit{spoof $k$-perfect number of the second kind} if $\sigma(n) (x^2+x+1) = knx$.
\end{enumerate}
Furthermore, we shall refer to $n$ as its \textit{root} and to $x$ as its \textit{spoof factor}.
\end{definition}
An advantage of this definition is that it is simple to implement as a computer program. Indeed, in order to find spoof $k$-perfect numbers, one merely has to run through positive integers $n$ and compute the quantity $q = \sigma(n)/kn$. If $q$ is of the form $x/(x+1)$ (for numbers of the first kind) or $x/(x^2+x+1)$ (for numbers of the second kind), then having $\gcd(x,x+1)=1$, respectively $\gcd(x,x^2+x+1)=1$, we can calculate the difference $\delta$ between its denominator and numerator,
$$
\delta = q_{den} - q_{num},
$$
and handle the following cases:
\begin{itemize}
    \item If $\delta = 1$, we have found a spoof $k$-perfect number of the first kind,
    \item If $\delta = q_{num}^2 + 1$, we have found a spoof $k$-perfect number of the second kind.
\end{itemize}
We implemented this algorithm in Wolfram Mathematica 11.1 and ran it up to $n=5\times 10^9$ with $2\leq k \leq 50$ for spoof odd $k$-perfect numbers of the first kind, and up to $n=10^8$ with $2\leq k \leq 100$ for the second kind. We found no such numbers of the first kind other than Descartes' example. On the other hand, we found several examples of the second kind. 

Table \ref{tabl-delta} shows the three odd spoof $k$-perfect numbers of the second kind that we found, for which $x$ is a prime that is also coprime to $n$.

\begin{table}[h!] 
\centering
\begin{tabular}{|c|c|c||c|} 
\hline
$s$  & $n$ & $x$ & $k$  \\ \hline
$77805$  & $11115$ & $7$ & $16$  \\ \hline
$92781$  & $1521$ & $61$ & $97$  \\ \hline
$8999757$  & $147537$ & $61$ & $98$  \\ \hline
\end{tabular}
\caption{Odd spoof $k$-perfect numbers $s=nx$ of the second kind} \label{tabl-delta}
\end{table}
Note that the spoof factor $61$ appears in two of our three examples, which is striking because it also appears in the spoof factor of Descartes' classical example.

Other odd examples exist, for which $x$ is not a prime number. For instance, $s=41089685$ (with $n=1173991$ and $x=35$) would be an odd $48$-perfect number if $x=35$ was the square of a prime.

Even spoof $k$-perfect numbers exist as well, and are much more abundant than their odd counterparts within our result set. For instance, $s=393120$ (with $n=10080, x=39$) would be a $4$-perfect number if only $39$ was prime and $s=1176725309760$ (with $n=862701840, x=1364$) would be a $5$-perfect number if only $1364$ was prime.

\section{On Robin's Inequality}

The spoof $k$-perfect examples we found and presented above are notable because the existence of multiperfect numbers is bounded by an inequality of Guy Robin \cite{Ro84}. Recall that we have 
$$
\sigma(n) < e^\gamma n \log \log n,
$$
for all positive integers $n > 5040$, where $\gamma$ denotes the Euler-Mascheroni constant, if and only if the Riemann Hypothesis is true. Thus, assuming the Riemann Hypothesis, we would expect a $k$-perfect number $n$ to appear only after
$$
n > e^{e^{k e^{-\gamma}}}.
$$
The spoof $k$-perfect numbers we found did not satisfy this inequality. For instance, we would expect a $98$-perfect number to appear after $e^{e^{98 e^{\gamma}}} \approx 10^{10^{23.53\ldots}}$, however our odd spoof $98$-perfect number,
$$
S = 8999757,
$$
seems to violate this bound. Other examples -- even and odd -- present such characteristics. For example, $s=941129280$ (with $n=470564640, x=2$), which is a spoof $7$-perfect number of the first kind, is found under the expected lower bound $1.3\times 10^{22}$ and $s=10805558400$ (with $n=2701389600, x=4$), a spoof $6$-perfect number also of the first kind, is found under the expected lower bound $4.1\times 10^{12}$, to name a few.

\section{Conclusion and Further Work}

In this paper we have introduced spoof $k$-perfect numbers as an extension to spoof perfect numbers and found several notable odd examples. The search for more of these numbers can be easily continued with more computing power using the same simple algorithm that we described above. Of course, one may extend our definition of such numbers even further by considering spoof factors with multiplicity greater than $2$. As an example, numbers $s=nx$ that would be perfect if only one of their prime factors were a \textit{cube} are such that
$$
\sigma(n)(x^3+x^2+x+1) = knx,
$$
and it would suffice to adapt our algorithm to search for positive integers $n$ satisfying
$$
\frac{\sigma(n)}{n} = \frac{kx}{x^3+x^2+x+1},
$$
for a positive integer $x$. In conclusion, we end this discussion with the following question, which is based on our experimental results.
\begin{question}
Does there exist an odd spoof $k$-perfect number of the first kind, for $k > 2$?
\end{question}

\vskip 20pt\noindent {\bf Acknowledgement.} The author is grateful to an anonymous referee for suggestions and corrections that improved the paper.


\begin{thebibliography}{1}\footnotesize

\bibitem{Pace spoof group} N. Andersen, S. Durham, M. Griffin,  J. Hales, P. Jenkins, R. Keck, H. Ko, G. Molnar, E. Moss, P. Nielsen, K. Niendorf, V. Tombs, M. Warnick, and D. Wu, Odd, spoof perfect factorizations,\textit{ J. Number Theory} \textbf{234}, (2022) 31-47.

\bibitem{Ro84} 
G. Robin, Grandes valeurs de la fonction somme des diviseurs et hypoth\`{e}se de Riemann, J. Math. Pures Appl. {\bf 63} (1984), 187--213.

\bibitem{Toth21}
L. T\'oth, On the Density of Spoof Odd Perfect Numbers, \textit{Comput. Methods Sci. Technol.} \textbf{27 (1)} (2021), 25--28.

\bibitem{Voight}
J. Voight, \textit{On the nonexistence of odd perfect numbers}, MASS selecta, Amer. Math. Soc., Providence, RI, 2003, pp. 293--300.

\end{thebibliography}
\end{document}